\documentclass [11pt,titlepage]{article}

\usepackage{amsmath}
\usepackage{amssymb}
\usepackage{mathrsfs}
\DeclareMathAlphabet{\mathpzc}{OT1}{pzc}{m}{it}
\newcommand{\BE}{\begin{equation}}
\newcommand{\EE}{\end{equation}}
\newcommand{\BEA}{\begin{eqnarray}}
\newcommand{\EEA}{\end{eqnarray}}
\newcommand{\BEAs}{\begin{eqnarray*}}
\newcommand{\EEAs}{\end{eqnarray*}}
\newcommand{\NN}{\nonumber}
\newcommand{\NI}{\noindent}
\newcommand{\DL}{\delta}
\newcommand{\EQ}{\equiv}
\newcommand{\LF}{\left}
\newcommand{\RI}{\right}
\newcommand{\LR}{\lfloor}
\newcommand{\RR}{\rfloor}
\newcommand{\ra}{\rightarrow}
\begin{document}
\begin{center}
{\large{\bf A unique method to evaluate the general integral $\int^{^{^\infty} }_{_{_0}}${\boldmath $\frac{ \sin ^{^a}px~\cos ^{^c}qx}{x^{^b}}~ dx$} }}\\
~\\
{\bf Joseph Amal Nathan}\\Reactor Physics Design Division, Bhabha Atomic Research Centre, Mumbai-400085, India.\\
{\bf email:}josephan@barc.gov.in
\end{center}

\NI
{\small {\bf Abstract:} All integrals available in literature and books, that are related to {\bf Sinc}$\left( \equiv \frac{\sin x}{x} \right)$ function {\bf \cite {GS}}, are special cases of the general form of the integral given in the title. The evaluation of the integral is divided into two cases (i) $a$ and $b$ of same parity, which is easier to evaluate and (ii) $a$ and $b$ of different parity, a difficult case. Amazingly and may be for the first time, a divergent integral is used in evaluating this difficult case with the help of a simple but a special combinatorial expression. The combinatorial identity is derived from the power reduction formula of the {\it sines} and {\it cosines}. The method adopted in this paper is unique and makes it relatively easy to evaluate this integral.

\NI
}

\NI
{\bf Introduction:} Some special cases of these integrals are classical. They have applications in signal processing, information theory and numerical analysis. In Mathematics they appear in the topic of Fourier Transforms(FT) and in Physics, in the study of diffraction of light by objects, called Fraunhofer diffraction(FD) {\bf \cite{JF}}. The FD is given by the FT of the object. For example the FD of a rectangular aperture (described by rectangular function), is proportional to the Sinc function. If we use the well known  Parseval's Theorem or Rayleigh's Identity (in Physics it is a consequence of Law of conservation of energy) we can evaluate, $\int^ \infty_0 \frac{ \sin^2x}{x^2} dx = \frac{\pi}{2}$, which can be integrated by parts to get $\int^ \infty_0 \frac{ \sin x}{x} dx = \frac{\pi}{2}$. If the object is a triangular aperture, then FD is proportional to Sinc$^2$, which can be used to evaluate, $\int^ \infty_0 \frac{ \sin ^4x}{x^4} dx = \frac{\pi}{3}$ and an integration by parts will give $\int^ \infty_0 \frac{ \sin ^3x}{x^3} dx = \frac{3 \pi}{8}$. 

\NI
The analytic form of the integral for any integer $a,b$ with $p=1,~c=0$, attributed to Prof. Wolstenholme, is discussed, in a classic book by Joseph Edwards {\bf \cite{JE}} and the expression for the same case can also be found in {\bf \cite{WR}}. The result for the case $a=b$ with $p=1,~c=0$ is derived in an innovative way in {\bf \cite{RB}}. The analytic expression for the special case $a=b=m$ with $c=p=1,~q=b$, is seen in {\bf \cite{AE}}. Since, analytic forms for only special cases are available in literature, the {{\bf motivation} is to find a form for a general integral.

\NI
For integers $a,b,c,p,q$ such that $a \geq b \geq 2$, $c \geq 0$, define,
$$
\mathcal I^{pq}_{abc} \EQ \LF[ \rule{0cm}{.35cm} \RI. \hbox{Sign}(p) \LF. \rule{0cm}{.35cm} \RI]^a \int^\infty_0 \frac{ \sin ^a|p|x~\cos ^c|q|x}{x^b} dx~~~~~~~\hbox{and}~~~~~~~\mathpzc d \! (x)^i_{ac} \EQ \frac{d^i}{dx^i} \LF[ \sin ^apx \cos ^cqx \RI].
$$
We know $\mathcal I^{10}_{ab0}$ converges. We have $\mathcal I^{p0}_{ab0}= p^{b-1} \mathcal I^{10}_{ab0}$ and $|\mathcal I^{pq}_{abc}| \leq |\mathcal I^{p0}_{ab0}|$. So $\mathcal I^{pq}_{abc}$ converges. Neglecting the sign of $p$ temporarily, there is no loss of generality if $p,q \geq 0$. Integrating by parts, the RHS of $\mathcal I^{pq}_{abc}$, $h$ times, such that $1 \leq h \leq (b-1)$, we get,
\BE \label{pii}
\mathcal I_{abc}^{pq}= - \sum_{i=0}^{h-1} \left[ \frac{(b-i-2)! ~\mathpzc d \! (x)^i_{ac}}{(b-1)!~ x^{b-i-1}} \right]^ \infty_0 + \frac{(b-h-1)!}{(b-1)!} \int ^ \infty _0 \mathpzc d \! (x)^h_{ac} ~\frac{dx}{x^{b-h}}.
\EE
The first part of eq(\ref{pii}) will have a term $\left[ \frac{\sin px}{x} \right]^{b-i-1} \!\! \left[ \sin px \right] ^{a-b+1}$, since $\LF[ \sin px \RI] ^{a-i}$ is a factor of $\mathpzc d \! (x)^i_{ac}$. This term is zero at $x=0$ and $x=\infty$, since $(b-i-1) \geq 1$, $(a-b+1) \geq 1$ and $| \sin px| \leq 1$. So eq(\ref{pii}) reduces to,
\BE \label{pil}
\mathcal I^{pq}_{abc} = \frac{(b-h-1)!}{(b-1)!} \int ^ \infty_0 \mathpzc d \! (x)^h_{ac} ~\frac{dx}{x^{b-h}}.
\EE
Let $s \EQ a ( \bmod~2)$ and $t \EQ c ( \bmod~2)$. Introducing integers $p,q$ in the trigonometric power-reduction formula, we get,
\BEA 
\label{soe} \sin ^a px =&& \!\!\!\!\!\!\!\!\!\! \frac{(1-s)}{2^a}{a \choose \frac{a}{2}} +\frac{2}{2^a} \sum_{i=0}^{ \LF \LR \!\! \frac{a-1}{2} \! \RI \RR} (-1)^{ \LF \LR \! \frac{a}{2} \! \RI \RR -i} {a \choose i} \LF[ s \sin (a-2i)px + (1-s) \cos (a-2i)px \RI], ~~~~~~~\\
\label{coe} \cos ^c qx =&& \!\!\!\!\!\!\!\!\!\! \frac{(1-t)}{2^c}{c \choose \frac{c}{2}} +\frac{2}{2^c} \sum_{i=0}^{\LF \LR \!\! \frac{c-1}{2} \! \RI \RR} {c \choose i} \cos (c-2i)qx, 
\EEA
where for any real $y,~ \LR \! y \! \RR$ is the largest integer not greater than $y$. Though eq(\ref{soe}) is not valid for $a=p=0$, since $a \geq 2$, we can use it. Let $a$ and $h$ are of opposite parity, using eq(\ref{soe}) and eq(\ref{coe}), along with the trigonometric product-to-sum identities, we get,
\BEA \label{aho}
\NN \mathpzc d \! (x)^h_{ac} \! =&& \!\!\!\!\!\!\!\!\!\! \frac{(-1)^{ \LF \LR \!\! \frac{h+1}{2} \! \RI \RR }}{2^{a+c-1}} \! \LF \{ \rule{0cm}{.7cm} \RI. \!\! (-1)^{ \LF \LR \! \frac{a}{2} \! \RI \RR } \!\! \LF[ \rule{0cm}{.5cm} \RI. \!\! (1-t) {c \choose \frac{c}{2}} \!\! \sum_{i=0}^{\LF \LR \!\! \frac{a-1}{2} \! \RI \RR } \!\! (-1)^i {a \choose i} [(a-2i)p]^h \sin (a-2i)px + \!\!\! \sum_{i=0}^{ \LF \LR \!\! \frac{a-1}{2} \! \RI \RR } \! \sum_{j=0}^{\LF \LR \!\! \frac{c-1}{2} \! \RI \RR } \! (-1)^i \\
\NN && \!\!\!\!\!\!\!\!\!\! {a \choose i} {c \choose j} \LF( \rule{0cm}{.36cm} \RI.  [(a-2i)p +(c-2j)q]^h \sin [(a-2i)p+(c-2j)q]x + [(a-2i)p-(c-2j)q]^h \\
&& \!\!\!\!\!\!\!\!\!\! \sin [(a-2i)p-(c-2j)q]x \LF. \rule{0cm}{.36cm} \RI) \LF. \rule{0cm}{.5cm} \RI] \! + \! (1-s) {a \choose \frac{a}{2}} \! \sum_{i=0}^{\LF \LR \!\! \frac{c-1}{2} \! \RI \RR } \! {c \choose i} [(c-2i)q]^h  \sin (c-2i)qx \LF. \rule{0cm}{.7cm} \RI \}.
\EEA
Similarly when $a$ and $h$ are of same parity, we have,
\BEA \label{ahs}
\NN \mathpzc d \! (x)^h_{ac} \! =&& \!\!\!\!\!\!\!\!\!\! \frac{(-1)^{ \LF \LR \! \frac{h}{2} \! \RI \RR }}{2^{a+c-1}} \! \LF \{ \rule{0cm}{.7cm} \RI. \!\! (-1)^{ \LF \LR \! \frac{a}{2} \! \RI \RR } \!\! \LF[ \rule{0cm}{.5cm} \RI. \!\! (1-t) {c \choose \frac{c}{2}} \!\! \sum_{i=0}^{\LF \LR \!\! \frac{a-1}{2} \! \RI \RR } \!\! (-1)^i {a \choose i} [(a-2i)p]^h \cos (a-2i)px + \!\!\! \sum_{i=0}^{ \LF \LR \!\! \frac{a-1}{2} \! \RI \RR } \! \sum_{j=0}^{\LF \LR \!\! \frac{c-1}{2} \! \RI \RR } \! (-1)^i \\
\NN && \!\!\!\!\!\!\!\!\!\! {a \choose i} {c \choose j} \LF( \rule{0cm}{.36cm} \RI.  [(a-2i)p +(c-2j)q]^h \cos [(a-2i)p+(c-2j)q]x + [(a-2i)p-(c-2j)q]^h \\
&& \!\!\!\!\!\!\!\! \cos [(a-2i)p-(c-2j)q]x \LF. \rule{0cm}{.36cm} \RI) \LF. \rule{0cm}{.5cm} \RI] \! + \! (1-s) {a \choose \frac{a}{2}} \! \sum_{i=0}^{\LF \LR \!\! \frac{c-1}{2} \! \RI \RR } \! {c \choose i} [(c-2i)q]^h  \cos (c-2i)qx \LF. \rule{0cm}{.7cm} \RI \}.
\EEA
Using eq(\ref{ahs}) we will now derive a general and a special combinatorial identity. Since $\mathpzc d \! (x)^h_{ac}$ contains the factor $\LF[ \sin px \RI] ^{a-h}$ with $(a-h) \geq 2$, from eq(\ref{ahs}) $\LF. \mathpzc d \! (x)^h_{ac} \RI|_{x= \pi} \! = 0$, so we get a general combinatorial identity,
\BEA \label{gce}
\NN && \!\!\!\!\!\!\!\!\!\!\!\!\! \mathpzc d \! (\pi)^h_{ac} = \! (-1) \! ^{\LF \LR \! \frac{a}{2} \! \RI \RR \! + ap}  \!\! \LF[ \rule{0cm}{.5cm} \RI. \!\!\! (1-t) \! {c \choose \frac{c}{2}} \!\!\! \sum_{i=0}^{\LF \LR \!\! \frac{a-1}{2} \! \RI \RR } \!\! (-1)^i \! {a \choose i} \! [(a-2i)p]^h \! + \! (-1)^{cq} \!\! \sum_{i=0}^{ \LF \LR \!\! \frac{a-1}{2} \! \RI \RR } \!\! \sum_{j=0}^{\LF \LR \!\! \frac{c-1}{2} \! \RI \RR } \! (-1)^i \! {a \choose i} \! {c \choose j} \!\! \LF( \rule{0cm}{.36cm} \RI. \!\! [(a-2i)p \\
&&~~~~~\! + (c-2j)q]^h \! + \! [(a-2i)p-(c-2j)q]^h \!\! \LF. \rule{0cm}{.36cm} \RI) \!\!\! \LF. \rule{0cm}{.5cm} \RI] \!\! + \! (1-s) \! (-1)^{cq} \! {a \choose \frac{a}{2}} \!\!\! \sum_{i=0}^{\LF \LR \!\! \frac{c-1}{2} \! \RI \RR } \!\!\! {c \choose i} \! [(c-2i)q]^h \!\! = \! 0.
\EEA
There will be four cases for eq(\ref{gce}). (i)$a$-even, $c$-even, we have, $(-1)^{ap}=(-1)^{cq}=1$. (ii)$a$-even, $c$-odd, the first term is zero, $(-1)^{ap}=1$ and the expression is independent of $(-1)^{cq}$. (iii)$a$-odd, $c$-even, the last term is zero, $(-1)^{cq}=1$ and the expression is independent of $(-1)^{ap}$. (iv)$a$-odd, $c$-odd, the first and last term is zero, and the expression is independent of $(-1)^{ap+cq}$. We later require the combinatorial expression with $h=b-1$. So when we put $h=b-1$ we see that, in eq(\ref{gce}), {\boldmath $a$} {\bf and} {\boldmath $b$} {\bf will be of opposite parity} and after accounting for the four cases discussed above we get the special combinatorial identity, which will be used below in {\it Case-1},
\BEA \label{sce} 
\NN \mathpzc d \! (\pi)^{b-1}_{ac} \!\!\!\!\!\!\!\!\!\! &&= \! (-1)^{ \LF \LR \! \frac{a}{2} \! \RI \RR } \!\!\! \LF[ \rule{0cm}{.5cm} \RI. \!\! (1-t) {c \choose \frac{c}{2}} \!\!\! \sum_{i=0}^{\LF \LR \!\! \frac{a-1}{2} \! \RI \RR } \! (-1)^i {a \choose i} [(a-2i)p]^{b-1} \!\! + \!\! \sum_{i=0}^{ \LF \LR \!\! \frac{a-1}{2} \! \RI \RR } \! \sum_{j=0}^{\LF \LR \!\! \frac{c-1}{2} \! \RI \RR } \! (-1)^i {a \choose i} \! {c \choose j} \!\! \LF( \rule{0cm}{.36cm} \RI. \!\! [(a-2i)p ~~~~~\\
&&~~ + \! (c-2j)q]^{b-1} \!\! + \! [(a-2i)p-(c-2j)q]^{b-1} \!\! \LF. \rule{0cm}{.36cm} \RI) \!\!\! \LF. \rule{0cm}{.5cm} \RI] \!\! + \! (1-s) \! {a \choose \frac{a}{2}} \!\!\! \sum_{i=0}^{\LF \LR \!\! \frac{c-1}{2} \! \RI \RR } \!\!\! {c \choose i} \! [(c-2i)q]^{b-1} \! = \! 0. 
\EEA
\NI
{\it Case-1:} {\boldmath $a > b$} {\it with opposite parity, and} {\boldmath $c \geq 0$}\\
\NI
We require the following, logarithmically divergent, integral. Let $l$ be an integer then,
\BEA \label{xsi}
\NN \int ^ \infty_0 \frac{\sin lx}{x^2} ~dx =&& \!\!\!\!\!\!\!\!\!\!\! \int^ \infty_0 \sin lx  \Bigg[ \! \int ^\infty_0 \!\! \Bigg( \!\! \int ^ \infty_v e^{-u x} ~ du \!\! \Bigg) dv \Bigg] dx = \int^ \infty_0 \!\! \Bigg( \!\! \int ^\infty_v \!\! \Bigg( \!\! \int ^ \infty_0 e^{-u x} \! \sin lx~ dx \!\! \Bigg) du \!\! \Bigg) dv ~~~~~\\
=&& \!\!\!\!\!\!\!\!\!\!\! \LF[ v \LF( \frac{\pi}{2}- \tan ^{-1} \frac{v}{l} \RI) + \frac{l}{2} \ln \LF( v^2+l^2 \RI) \RI]^ \infty _0 = l \LF( 1+ \lim_{v \ra \infty} \ln v - \ln \LF[ |l| + \DL_{|l|0} \RI] \RI),
\EEA 
where $\DL_{|l|0}$ is added for the case $l=0$.

\NI
Substituting eq(\ref{aho}) with $h=b-2$ in eq(\ref{pil}) we get,
\BEAs
\NN \mathcal I^{pq}_{abc} \! =&& \!\!\!\!\!\!\!\!\!\! \frac{(-1)^{ \LF \LR \!\! \frac{b-1}{2} \! \RI \RR }}{2^{a+c-1}(b-1)!} \! \LF \{ \rule{0cm}{.7cm} \RI. \!\! (-1)^{ \LF \LR \! \frac{a}{2} \! \RI \RR } \!\! \LF[ \rule{0cm}{.5cm} \RI. \!\! (1-t) {c \choose \frac{c}{2}} \!\! \sum_{i=0}^{\LF \LR \!\! \frac{a-1}{2} \! \RI \RR } \!\! (-1)^i {a \choose i} [(a-2i)p]^{b-2} \int ^ \infty_0 \frac{\sin (a-2i)px}{x^2} dx ~~~~~~~~\\
\NN && \!\!\!\!\!\!\!\!\! + \sum_{i=0}^{ \LF \LR \!\! \frac{a-1}{2} \! \RI \RR } \! \sum_{j=0}^{\LF \LR \!\! \frac{c-1}{2} \! \RI \RR } \! (-1)^i {a \choose i} {c \choose j} \LF( \rule{0cm}{.36cm} \RI.  [(a-2i)p +(c-2j)q]^{b-2} \int ^ \infty_0 \frac{\sin [(a-2i)p+(c-2j)q]x}{x^2} dx \\
\NN && \!\!\!\!\!\!\!\!\! + [(a-2i)p-(c-2j)q]^{b-2} \int ^ \infty_0 \frac{\sin [(a-2i)p-(c-2j)q]x}{x^2} dx \LF. \rule{0cm}{.36cm} \RI) \LF. \rule{0cm}{.5cm} \RI] \!+\! (1-s) {a \choose \frac{a}{2}} \! \sum_{i=0}^{\LF \LR \!\! \frac{c-1}{2} \! \RI \RR }\! {c \choose i} \\
\NN && \!\!\!\!\!\!   [(c-2i)q]^{b-2}  \int ^ \infty_0 \frac{\sin (c-2i)qx}{x^2} dx \LF. \rule{0cm}{.7cm} \RI \}.
\EEAs
Since $p,q$ are positive, in the above equation the coefficients of the {\it Sine} functions are positive, except $(a-2i)p-(c-2j)q$, which {\bf may} take zero or negative values. The positive coefficients are zero only when $p$ and/or $q$ are zero. Substituting eq(\ref{xsi}) in the above equation we get,
\BEAs
\NN \mathcal I^{pq}_{abc} \! =&& \!\!\!\!\!\!\!\!\!\! \frac{(-1)^{ \LF \LR \!\! \frac{b-1}{2} \! \RI \RR }}{2^{a+c-1}(b-1)!} \! \LF \{ \rule{0cm}{.7cm} \RI. \!\! (-1)^{ \LF \LR \! \frac{a}{2} \! \RI \RR } \!\! \LF[ \rule{0cm}{.5cm} \RI. \!\! (1-t) {c \choose \frac{c}{2}} \!\! \sum_{i=0}^{\LF \LR \!\! \frac{a-1}{2} \! \RI \RR } \!\! (-1)^i {a \choose i} [(a-2i)p]^{b-1} \LF( \rule{0cm}{.36cm} \RI. 1+ \lim_{v \ra \infty} \ln v ~~~~~~~~
\EEAs
\BEAs
\NN && \!\!\!\!\! - \ln \LF[(a-2i)p + \DL_{p0} \RI] \LF. \rule{0cm}{.36cm} \RI) + \sum_{i=0}^{ \LF \LR \!\! \frac{a-1}{2} \! \RI \RR } \! \sum_{j=0}^{\LF \LR \!\! \frac{c-1}{2} \! \RI \RR } \! (-1)^i {a \choose i} {c \choose j} \LF[ \rule{0cm}{.36cm} \RI.  [(a-2i)p +(c-2j)q]^{b-1} \\
\NN && \!\!\!\!\! \LF( 1+ \lim_{v \ra \infty} \ln v - \ln \LF[ (a-2i)p +(c-2j)q + \DL_{p0} \DL_{q0} \RI] \RI)+ [(a-2i)p-(c-2j)q]^{b-1} \\
\NN && \!\!\!\!\! \LF( 1+ \lim_{v \ra \infty} \ln v - \ln \LF[ |(a-2i)p -(c-2j)q| + \DL_{|(a-2i)p -(c-2j)q|0} \RI] \RI) \LF. \rule{0cm}{.36cm} \RI] \LF. \rule{0cm}{.5cm} \RI] \!+\! (1-s) {a \choose \frac{a}{2}} \\
\NN && \!\!\!\! \sum_{i=0}^{\LF \LR \!\! \frac{c-1}{2} \! \RI \RR }\! {c \choose i} [(c-2i)q]^{b-1}  \LF( 1+ \lim_{v \ra \infty} \ln v - \ln \LF[ (c-2i)q + \DL_{q0} \RI] \RI) \LF. \rule{0cm}{.7cm} \RI \},
\EEAs
which can be rewritten using eq(\ref{sce}) as 
\BEAs
\NN \mathcal I^{pq}_{abc} \! =&& \!\!\!\!\!\!\!\!\!\! \frac{(-1)^{ \LF \LR \!\! \frac{b-1}{2} \! \RI \RR }}{2^{a+c-1}(b-1)!} \! \LF \{ \rule{0cm}{.7cm} \!\!\!\! \RI. \LF( \!  1+ \lim_{v \ra \infty} \ln v \! \RI) \mathpzc d \! (\pi)^{b-1}_{ac} - (-1)^{ \LF \LR \! \frac{a}{2} \! \RI \RR } \!\! \LF[ \rule{0cm}{.5cm} \RI. \!\! (1-t) {c \choose \frac{c}{2}} \!\! \sum_{i=0}^{\LF \LR \!\! \frac{a-1}{2} \! \RI \RR } \!\! (-1)^i {a \choose i} [(a-2i)p]^{b-1}  ~~~~~~~~\\
\NN && \!\!\!\!\!\!\!\!\! \ln \LF[(a-2i)p + \DL_{p0} \RI] + \sum_{i=0}^{ \LF \LR \!\! \frac{a-1}{2} \! \RI \RR } \! \sum_{j=0}^{\LF \LR \!\! \frac{c-1}{2} \! \RI \RR } \! (-1)^i {a \choose i} {c \choose j} \LF( \rule{0cm}{.36cm} \RI.  [(a-2i)p +(c-2j)q]^{b-1} \ln \LF[ (a-2i)p \RI. \\
\NN && \!\!\!\!\!\!\!\!\! +(c-2j)q + \LF. \DL_{p0} \DL_{q0} \RI] + [(a-2i)p-(c-2j)q]^{b-1} \ln \LF[ |(a-2i)p -(c-2j)q| \RI. \\
\NN && \!\!\!\!\!\!\!\!\! \LF. + \DL_{|(a-2i)p -(c-2j)q|0} \RI] \LF. \rule{0cm}{.36cm} \RI) \LF. \rule{0cm}{.5cm} \RI] \!-\! (1-s) {a \choose \frac{a}{2}} \sum_{i=0}^{\LF \LR \!\! \frac{c-1}{2} \! \RI \RR }\! {c \choose i} [(c-2i)q]^{b-1}  \ln \LF[ (c-2i)q + \DL_{q0} \RI] \LF. \rule{0cm}{.7cm} \RI \}.
\EEAs
Again from eq(\ref{sce}) $\mathpzc d \! (\pi)^{b-1}_{ac}=0$, so for {\it Case-1:}, with the sign of $p$, the above equation becomes,
\BEA \label{ca1}  
\NN \mathcal I^{pq}_{abc} =&& \!\!\!\!\!\!\!\!\! \frac{\LF[ \hbox{Sign}(p) \RI]^a(-1)^{ \LF \LR \!\! \frac{b+1}{2} \! \RI \RR  }}{2^{a+c-1}(b-1)!} \LF \{ \rule{0cm}{.7cm} \RI. (1-t) (-1)^{ \LF \LR \! \frac{a}{2} \! \RI \RR } {c \choose \frac{c}{2}} \sum_{i=0}^{\LF \LR \!\! \frac{a-1}{2} \! \RI \RR } (-1)^i {a \choose i} [(a-2i)p]^{b-1} \ln [(a-2i)p+\DL_{p0}] \\
\NN && \!\!\!\!\!\!\!\!\! + (-1)^{ \LF \LR \! \frac{a}{2} \! \RI \RR } \sum_{i=0}^{ \LF \LR \!\! \frac{a-1}{2} \! \RI \RR } \sum_{j=0}^{\LF \LR \!\! \frac{c-1}{2} \! \RI \RR } (-1)^i {a \choose i} {c \choose j} \LF( \rule{0cm}{.4cm} \RI. [(a-2i)p+(c-2j)q]^{b-1} \ln[(a-2i)p\\
\NN && \!\!\!\!\!\!\!\!\! +(c-2j)q+ \DL_{p0} \DL_{q0}] + [(a-2i)p-(c-2j)q]^{b-1} \ln \LF[ \RI. |(a-2i)p-(c-2j)q|\\
&& \!\!\!\!\!\!\!\!\! + \DL_{|(a-2i)p-(c-2j)q|0} \LF. \RI] \LF. \rule{0cm}{.4cm} \RI) \!+\! (1-s) {a \choose \frac{a}{2}} \! \sum_{i=0}^{\LF \LR \!\! \frac{c-1}{2} \! \RI \RR } \! {c \choose i} [(c-2i)q]^{b-1} \! \ln [(c-2i)q+ \DL_{q0}] \LF. \rule{0cm}{.7cm} \RI \}.  
\EEA
\NI
{\it Case-2:} {\boldmath $a \geq b$} {\it with same parity, and} {\boldmath $c \geq 0$}\\
\NI
We will evaluate the following integral. Let $k$ be an integer,
\BEA \label{s/x}
\NN \int ^ \infty_0 \frac{\sin kx}{x} ~dx &=& \int ^ \infty_0 \sin kx \LF ( \int ^ \infty_0 e^{-ux} ~du \RI ) ~dx = \int ^ \infty_0 \LF ( \int ^ \infty_0 e^{-ux} \sin kx ~dx \RI ) ~du,\\
&=& \LF[ \tan^{-1} \LF( \frac{u}{k} \RI) \RI ]_0^\infty = \hbox{Sign}(k) (1-\DL_{|k|0})\frac{\pi}{2},
\EEA
where $\DL_{|k|0}$ is introduced for the case $k=0$. 

\NI
With $h=b-1$, substituting eq(\ref{aho}) in eq(\ref{pil}) we get,
\BEAs
\NN \mathcal I^{pq}_{abc} \! =&& \!\!\!\!\!\!\!\!\!\! \frac{(-1)^{ \LF \LR \! \frac{b}{2} \! \RI \RR }}{2^{a+c-1}(b-1)!} \! \LF \{ \rule{0cm}{.7cm} \RI. \!\! (-1)^{ \LF \LR \! \frac{a}{2} \! \RI \RR } \!\! \LF[ \rule{0cm}{.5cm} \RI. \!\! (1-t) {c \choose \frac{c}{2}} \!\! \sum_{i=0}^{\LF \LR \!\! \frac{a-1}{2} \! \RI \RR } \!\! (-1)^i {a \choose i} [(a-2i)p]^{b-1} \int ^ \infty_0 \frac{\sin (a-2i)px}{x} dx ~~~~~~~~\\
\NN && \!\!\!\!\!\!\!\!\!\! + \sum_{i=0}^{ \LF \LR \!\! \frac{a-1}{2} \! \RI \RR } \! \sum_{j=0}^{\LF \LR \!\! \frac{c-1}{2} \! \RI \RR } \! (-1)^i {a \choose i} \! {c \choose j} \! \LF( \rule{0cm}{.36cm} \RI.  [(a-2i)p +(c-2j)q]^{b-1} \! \int ^ \infty_0 \frac{\sin [(a-2i)p+(c-2j)q]x}{x} dx \\
\NN && \!\!\!\!\!\!\!\!\!\! + [(a-2i)p-(c-2j)q]^{b-1} \int ^ \infty_0 \frac{\sin [(a-2i)p-(a-2j)q]x}{x} dx \LF. \rule{0cm}{.36cm} \RI) \LF. \rule{0cm}{.5cm} \RI] \! \\
\NN && \!\!\!\!\!\!\!\!\!\! + (1-s) {a \choose \frac{a}{2}} \! \sum_{i=0}^{\LF \LR \!\! \frac{c-1}{2} \! \RI \RR } \! {c \choose i} [(c-2i)q]^{b-1}  \int ^ \infty_0 \frac{\sin (c-2i)qx}{x} dx \LF. \rule{0cm}{.7cm} \RI \},
\EEAs
which after using eq(\ref{s/x}), with the sign of $p$, becomes,
\BEA \label{ca2}
\NN \mathcal I^{pq}_{abc} =&& \!\!\!\!\!\!\!\!\!\! \frac{\LF[ \hbox{Sign}(p) \RI]^a(-1)^{ \LF \LR \! \frac{b}{2} \! \RI \RR } \pi}{2^{a+c}(b-1)!} \LF \{ \rule{0cm}{.7cm} \RI. (1-t) (-1)^{ \LF \LR \! \frac{a}{2} \! \RI \RR } {c \choose \frac{c}{2}} \sum_{i=0}^{\LF \LR \!\! \frac{a-1}{2} \! \RI \RR } (-1)^i {a \choose i} [(a-2i)p]^{b-1} + (-1)^{ \LF \LR \! \frac{a}{2} \! \RI \RR } \sum_{i=0}^{ \LF \LR \!\! \frac{a-1}{2} \! \RI \RR } \\
\NN && \!\!\!\!\!\!\!\!\!\! \sum_{j=0}^{\LF \LR \!\! \frac{c-1}{2} \! \RI \RR } (-1)^i {a \choose i} {c \choose j} \LF( \rule{0cm}{.36cm} \RI. [(a-2i)p+(c-2j)q]^{b-1} + \hbox{Sign} \LF( \rule{0cm}{0.35cm} \RI. (a-2i)p-(c-2j)q \LF. \rule{0cm}{0.35cm} \RI) \\
&& \!\!\!\!\!\!\!\!\! [(a-2i)p-(c-2j)q]^{b-1} \LF. \rule{0cm}{.36cm} \RI) + (1-s) {a \choose \frac{a}{2}} \sum_{j=0}^{\LF \LR \!\! \frac{c-1}{2} \! \RI \RR } {c \choose j} [(c-2j)q]^{b-1} \LF. \rule{0cm}{.7cm} \RI \}.
\EEA
Eq(\ref{ca1}) and eq(\ref{ca2}) gives the complete analytic form of $\mathcal I^{pq}_{abc}$.

\NI
{\bf Acknowledgment:} I thank M.A. Prasad for many fruitful discussions and suggestions.


\begin{thebibliography}{99}
\bibitem{RB}R. Butler, {\it On the evaluation of the integral} $ \int ^ \infty_0 \sin^m t/x^m dt$, Amer. Math. Monthly, v67, 1960, pp. 566-569.
\bibitem{AE}A. Erdélyi, W. Magnus, P. Oberhettinger and F. G. Tkicomi, {\it Tables of Integral Transforms Vol-I}, McGraw-Hill, New York, 1954, p. 20.
\bibitem{GS}W. B. Gearhart and H. S. Schultz, {\it The function} $\frac{\sin x}{x}$, College Math J., {\bf 21}, 1991, pp. 90-99.
\bibitem{JF}J. F. James, {\it A Student's Guide to Fourier Transforms}, Cambridge University Press, 2011, pp. 40-55.
\bibitem{JE}Joseph Edwards, {\it A treatise on the Integral Calculus Vol-II}, Macmillan and Co. Ltd, 1922, pp. 188-213.
\bibitem{WR}Weisstein, Eric W., {\it Sinc Function}, From MathWorld-A Wolfram Web Resource. (http://mathworld.wolfram.com/SincFunction.html)
\end{thebibliography}
\end{document}